\documentclass[11pt]{amsart}
\usepackage{palatino}
\usepackage{amsfonts}
\usepackage{amssymb}
\usepackage{amscd}

\theoremstyle{plain}

\newtheorem{theorem}{Theorem}[section]
\newtheorem{proposition}[theorem]{Proposition}

\theoremstyle{definition}

\newtheorem{remark}[theorem]{Remark}

\newtheorem{conjecture/question}[theorem]{Conjecture/Question}

\newtheorem{remark/definition}[theorem]{Remark/Definition}
\newtheorem{terminology/notation}[theorem]{Terminology/Notation}

\newcommand{\marginlabel}[1]%
  {\mbox{}\marginpar{\raggedleft\hspace{0pt}\bfseries\sf#1}}

\def\PP{{\textbf P}}
\def\OO{\mathcal{O}}

\def\H{\mathcal{H}}

\def\cM{\mathcal{M}}

\def\cU{\mathcal{U}}

\def\H{\mathcal{H}^{\mathrm{tr}}}
\def\Pic0{{\rm Pic}^0(X)}

\def\mm{\overline{\mathcal{M}}}

\def\dd{\overline{\mathcal{D}}}
\def\tr{\overline{\mathfrak{TR}}}
\def\ww{\overline{\mathcal{W}}}
\def\hh{\overline{\mathcal{H}}^{\mathrm{tr}}}
\def\aa{\overline{\mathcal{A}}}
\begin{document}

\title{\bf The Fermat cubic and special Hurwitz loci in $\mm_g$}

\author[G. Farkas]{Gavril Farkas}
\address{Humboldt Universit\"at zu Berlin, Institut f\"ur Mathematik,
10099 Berlin} \email{{\tt farkas@math.hu-berlin.de}}
\thanks{Research  partially supported by an Alfred P. Sloan Fellowship and the NSF Grant DMS-0500747
}

\maketitle

\pagestyle{myheadings}
\theoremstyle{remark}
\newtheorem*{nctheorem}{Theorem}

{\small{Abstract: We compute the class of the compactified Hurwitz
divisor $\tr_d$ in $\mm_{2d-3}$ consisting of curves of genus
$g=2d_3$ having a pencil $\mathfrak g^1_d$ with two unspecified
triple ramification points. This is the first explicit example of a
geometric divisor on $\mm_g$ which is not pulled-back form the
moduli space of pseudo-stable curves. We show  that the intersection
of $\tr_d$ with the boundary divisor $\Delta_1$ in $\mm_g$ picks-up
the locus of Fermat cubic tails.}}

\section{Introduction}

Hurwitz loci have played a basic role in the study of the moduli
space of curves at least since 1872 when Clebsch, and later Hurwitz,
proved that $\cM_g$ is irreducible by showing that a certain Hurwitz
space  parameterizing coverings of $\PP^1$ is connected (see
\cite{Hu}, or \cite{Fu2} for a modern proof). Hurwitz cycles on
$\mm_g$ are essential in the work of Harris and Mumford \cite{HM} on
the Kodaira dimension of $\mm_g$ and
 are expected to govern the length of minimal affine stratifications of
$\cM_g$. Faber and Pandharipande have proved that the class of any
Hurwitz cycle on $\mm_{g, n}$ is tautological (cf. \cite{FP}).  Very
few explicit formulas for the classes of such cycles are known.

We define a \emph{Hurwitz divisor in} $\mm_g$  \emph{with} $n$
\emph{degrees of freedom} as follows: We fix integers $k_1, \ldots,
k_n\geq 3$ and positive integers $d, g$ such that
$$k_1+k_2+\cdots+k_n=2d-g+n-1.$$ Then $\mathcal{H}_{g: \ k_1,
\ldots, k_n}$ is the locus of curves $[C]\in \cM_g$ having a degree
$d$ morphism $f:C\rightarrow \PP^1$ together with $n$ distinct
points $p_1, \ldots, p_n \in C$ such that $\mbox{mult}_{p_i}(f)\geq
k_i$ for $i=1, \ldots, n$. When $n=0$ and $g=2d-1$, we recover the
Brill-Noether divisor of $d$-gonal curves  studied extensively in
\cite{HM}. For $n=1$ we obtain Harris' divisor $\mathcal{H}_{g:\ k}$
of curves having a linear series $C\stackrel{d:1}\rightarrow \PP^1$
with a $k=(2d-g+1)$-fold point, cf. \cite{H}. If $n=1$ and $d=g-1$
then $\mathcal{H}_{g:\ g-1}$ specializes to  S. Diaz's divisor of
curves $[C]\in \cM_g$ having an exceptional Weierstrass point $p\in
C$ with $h^0(C, \OO_C((g-1)p))\geq 1$ (cf. \cite{Di}).

Since $\mathcal{H}_{g: k_1, \ldots, k_n}$ is the push-forward of a
cycle of codimension $n+1$ in $\cM_{g, n}$, as $n$ increases the
problem of calculating the class of $\overline{\mathcal{H}}_{g: k_1,
\ldots, k_n}$ becomes more and more difficult. In this paper we
carry out the first study of a Hurwitz locus having at least $2$
degrees of freedom, and we treat the simplest non-trivial case, when
$n=2$, $k_1=k_2=3$ and $g=2d-3$. Our main result is the calculation
of the class of $\tr_d:=\overline{\mathcal{H}}_{2d-3: \ 3, 3}$. As
usual we denote by $\lambda\in \mbox{Pic}(\mm_g)$ the Hodge class
and by $\delta_0, \ldots, \delta_{[g/2]}\in \mbox{Pic}(\mm_g)$ the
codimension $1$ classes on the moduli stack corresponding to the
boundary divisors of $\mm_g$:

\begin{theorem}\label{classoftr}
We fix $d\geq 3$ and denote by $\mathfrak{TR}_d$ the locus of curves
$[C]\in \cM_{2d-3}$ having a covering $C\stackrel{d:1}\rightarrow
\PP^1$ with  two unspecified triple ramification points. Then
$\mathfrak{TR}_d$ is an effective divisor on $\cM_{2d-3}$ and the
class of its compactification $\tr_d$ inside $\mm_{2d-3}$ is given
by the formula:
$$\tr_d \equiv 2\frac{(2d-6)!}{d!\ (d-3)!}\bigl(a\ \lambda-b_0\ \
\delta_0-b_1\ \delta_1-\cdots-b_{d-2}\ \delta_{d-2}\bigr)\in
\mathrm{Pic}(\mm_{2d-3}),$$ where
$$a=24(36d^4-36d^3-640d^2+1885-1475), \ \
b_0=144d^4-528d^3-298d^2+3049d-2940,$$
$$ \mbox{and }\ \  b_i=12i
(2d-3-i)(36d^3-156d^2+180d-5), \ \mbox{ for } \ 1\leq i\leq d-2.$$
\end{theorem}

The divisor $\tr_d$ is also the first example of a geometric divisor
in $\mm_g$ which is not a pull-back of an effective divisor from the
space $\mm_g^{\mathrm{ps}}$ of pseudo-stable curves. Precisely, if
we denote by $R\subset \mm_g$ the extremal ray obtained by attaching
to a fixed pointed curve $[C, q]$ of genus $g-1$ a pencil of plane
cubics, then $R\cdot \lambda=1, R\cdot \delta_0=12$, $R\cdot
\delta_1=-1$ and $R\cdot \delta_{\alpha}=0$ for $\alpha \geq 2$. If
$\delta:=\delta_0+\cdots+\delta_{[g/2]}\in \mbox{Pic}(\mm_g)$ is the
total boundary, there exists a divisorial contraction of the
extremal ray $R\subset \Delta_1 \subset \mm_g$ induced by the base
point free linear system $|11\lambda-\delta|$ on $\mm_g$,
$$f:\mm_g \rightarrow \mm_g^{\mathrm{ps}}.$$
The image is isomorphic to the moduli space of pseudo-stable curves
as defined by D. Schubert in \cite{S}. A curve is
\emph{pseudo-stable} if it has only nodes and cusps as
singularities, and each component of genus $1$ (resp. $0$)
intersects the curve in at least $2$ (resp. $3$ points). The
contraction $f$ is the first step in carrying out the minimal model
program for $\mm_g$, see \cite{HH}. One has an inclusion
$f^*(\mbox{Eff}(\mm_g^{\mathrm{ps}}))\subset \mbox{Eff} (\mm_g)$.
All the geometric divisors on $\mm_g$ whose class has been computed
(e.g. Brill-Noether or Gieseker-Petri divisors \cite{EH}, Koszul
divisors \cite{Fa}, \cite{Fa2}, or loci of curves with an abnormal
Weierstrass point \cite{Di}), lie in the subcone
$f^*(\mbox{Eff}(\mm_g^{\mathrm{ps}}))$. The divisor $\tr_d$ behaves
quite differently: If $i:\Delta_1=\mm_{1, 1}\times \mm_{g-1,
1}\hookrightarrow \mm_g$ denotes the inclusion, then we have the
relation
$$i^*(\tr_{d})= \alpha\cdot \{j=0 \}
\times \mm_{g-1, 1}+\mm_{1, 1}\times D=\alpha\cdot
\bigl\{\mbox{Fermat cubic}\bigr\}\times \mm_{g-1, 1}+\mm_{1,
1}\times D,$$ where $\alpha :=\frac{3 (2d-4)!}{d!\ (d-3)!}$ and
$D\subset \mm_{g-1, 1}$ is an explicitly described effective
divisor. Hence when restricted to the boundary divisor
$\Delta_1\subset \mm_g$ of elliptic tails, $\tr_d$ picks-up the
locus of \emph{Fermat cubic tails}!

The rich geometry of $\tr_d$ can also be seen at the level of genus
$2$ curves. We denote by $\chi:\mm_{2, 1}\rightarrow \mm_{2d-3}$ be
the map obtained by attaching a fixed tail $[B, q]$ of genus $2d-5$
at the marked point of every curve of genus $2$. Then the pull-back
under $\chi$ of every known geometric divisor on $\mm_{2, 1}$ is a
multiple of the Weierstrass divisor $\ww$ of $\mm_{2, 1}$ (cf.
\cite{HM}, \cite{EH}, \cite{Fa}). In contrast, for $\tr_d$ we have
the following picture:

\begin{theorem}\label{genus2} If $\chi:\mm_{2, 1}\rightarrow \mm_g$ is as above,
we have the following relation in $\mathrm{Pic}(\mm_{2, 1})$:
$$\chi^*(\tr_d)=N_1(d)\cdot \ww +e(d, 2d-5)\cdot
\overline{\mathcal{D}}_1+a(d-1, 2d-5)\cdot
\overline{\mathcal{D}}_2+a(d, 2d-5)\cdot \overline{\mathcal{D}}_3,$$
$$\mbox{ where } \ \mathcal{W}:=\{[C, p]\in \cM_{2, 1}: p\in C
\mbox{ is a Weierstrass point}\},$$
$$\mathcal{D}_1:=\{[C, p]\in \cM_{2, 1}: \exists x\in C-\{p \}
\mbox{ such that } 3x\equiv 3p\},$$
$$\mathcal{D}_2:=\{[C, p]\in \cM_{2, 1}:\exists l\in G^1_3(C), x\neq y\in C-\{p \}\ \mbox{
 with } a_1^l(x)\geq 3, a_1^l(y)\geq 3,\  a_1^l(p)\geq 2\},$$ and
$$\mathcal{D}_3:=\{[C, p]\in \cM_{2, 1}:\exists l\in G^1_4(C), x\neq
y\in C-\{p\} \mbox{ with } a_1^l(p)\geq 4, \ a_1^l(x)\geq 3,\
a_1^l(y)\geq 3 \}.$$
\end{theorem}

The constants $N_1(d), e(d, 2d-5), a(d, 2d-5), a(d-1, 2d-5)$
appearing in the statement are explicitly known and defined in
Proposition \ref{harris}. We used the notation
$a_1^l(p):=\mbox{mult}_p(l)$, for the multiplicity of a pencil $l\in
G^1_d(C)$ at a point $p\in C$. The classes of the divisors $\dd_1,
\dd_2, \dd_3$ on $\mm_{2, 1}$ are determined as well (The class of
$\ww$ is of course well-known, see \cite{EH}):
\begin{theorem}\label{3div}
One has the following formulas expressed in the basis $\{\psi,
\lambda, \delta_0\}$ of  $\mathrm{Pic}(\mm_{2, 1})$:
$$\dd_1\equiv 80\psi+10\delta_0-120\lambda,\ \ \mbox{ }\ \ \dd_2\equiv 160\psi+17\delta_0-200\lambda,$$
$$\ \mbox{ and }\ \ \
\dd_3\equiv 640\psi+72\delta_0-860\lambda.$$
\end{theorem}

\noindent {\bf{Acknowledgment:}} I have benefitted from discussions
with R. Pandharipande (5 years ago!) on counting admissible
coverings.

\section{Admissible coverings with two triple points}

We begin by recalling a few facts about admissible coverings in the
context of points of triple ramification. Let $\H_d$ be the Hurwitz
space parameterizing degree $d$ maps $[f:C\rightarrow \PP^1, q_1,
q_2; p_1, \ldots, p_{6d-12}]$, where $[C]\in \cM_{2d-3}$, \ $q_1,
q_2, p_1, \ldots, p_{6d-12}$ are distinct points on $\PP^1$ and $f$
has one point of triple ramification over each of $q_1$ and $q_2$
and one point of simple ramification over $p_i$ for $1\leq i\leq
6d-12$. We denote by $\hh_d$ the compactification of the Hurwitz
space by means of Harris-Mumford admissible coverings (cf.
\cite{HM}, \cite{ACV} and \cite{Di} Section 5; see also \cite{BR}
for a survey on Hurwitz schemes and their compactifications). Thus
$\hh_d$ is the parameter space of degree $d$ maps
$$[f:X\stackrel{d:1}\longrightarrow R, q_1, q_2; p_1, \ldots,
p_{6d-12}],$$  where $[R, q_1, q_2; p_1, \ldots, p_{6d-12}]$ is a
nodal rational curve, $X$ is a nodal curve of genus $2d-3$ and $f$
is a finite map which satisfies the following conditions:

\noindent $\bullet$ $f^{-1}(R_{\mathrm{reg}})=X_{\mathrm{reg}}$ and
$f^{-1}(R_{\mathrm{sing}})=X_{\mathrm{sing}}$.

\noindent $\bullet$ $f$ has a point of triple ramification over each
of $q_1$ and $q_2$ and simple ramification over $p_1, \ldots,
p_{6d-12}$. Moreover $f$ is \'etale over each point in
$R_{\mathrm{reg}}-\{q_1, q_2, p_1, \ldots, p_{6d-12}\}$.

\noindent $\bullet$ If $x\in X_{\mathrm{sing}}$ and $x\in X_1\cap
X_2$ where $X_1$ and $X_2$ are irreducible components of $X$, then
$f(X_1)$ and $f(X_2)$ are distinct components of $R$ and
$$\mathrm{mult}_{x} \{f_{| X_1}:X_1\rightarrow f(X_1)\}=\mathrm{mult}_{x} \{f_{|
X_2}:
X_2\rightarrow f(X_2)\}.$$

The group $\mathfrak{S}_2\times \mathfrak{S}_{6d-12}$ acts on
$\hh_d$ by permuting the triple and the ordinary ramification points
of $f$ respectively and we denote by $\mathfrak
H_d:=\hh_d/\mathfrak{S}_2\times \mathfrak{S}_{6d-12}$ for the
quotient. There exists a stabilization morphism $\sigma:
\mathfrak{H}_d\rightarrow \mm_g$ as well as a  finite map
$\beta:\mathfrak{H}_d\rightarrow \mm_{0, 6d-10}$. The description of
the local rings of $\hh_d$ can be found in \cite{HM} pg. 61-62 or
\cite{BR} and will be used in the paper. In particular, the scheme
$\hh_d$ is smooth at points $[f:X\rightarrow R, q_1, q_2; p_1,
\ldots, p_{6d-12}]$ with the property that there are no
automorphisms $\phi:X\rightarrow X$ with $f\circ \phi=f$.

\subsection{The enumerative geometry of pencils on the general curve}
We shall determine the intersection multiplicities of $\tr_d$ with
standard test curves in $\mm_g$. For this we need a variety of
enumerative results concerning pencils on pointed curves which will
be used throughout the paper. For a point $p\in C$ and a linear
series $l\in G^r_d(C)$, we denote by $$a^l(p):
\bigl(0<a_0^l(p)<a_1^l(p)<\ldots < a_r^l(p)\leq d\bigr)$$ the
\emph{vanishing sequence } of $l$ at $p$. If $l\in G^1_d(C)$, we say
that $p\in C$ is an \emph{$n$-fold point} if $l(-np)\neq \emptyset$.
We first recall the results from \cite{HM} Theorem A and \cite{H}
Theorem 2.1.

\begin{proposition}\label{harris}
Let us fix a general curve $[C, p]\in \cM_{g, 1}$  and an integer
$d\geq 2d-g-1\geq 0$. \newline \noindent $\bullet$ The number of
pencils $L\in W^1_d(C)$ satisfying $h^0(L\otimes \OO_C(-(2d-g-1)
p))\geq 1$ equals
$$a(d, g):=(2d-g-1)\frac{g!}{d!\ (g-d+1)!}.$$
\newline
\noindent $\bullet$ The number of pairs $(L, x)\in W^1_d(C)\times C$
satisfying $h^0(L\otimes \OO_C(-(2d-g) x))\geq 2$ equals
$$b(d, g):=(2d-g-1)(2d-g)(2d-g+1)\frac{g!}{d!\ (g-d)!}.$$
\newline
\noindent $\bullet$ Fix integers $\alpha, \beta\geq 1$ such that
$\alpha+\beta=2d-g$. The number of pairs $(L, x)\in W^1_d(C)\times
C$ satisfying $h^0(L\otimes \OO_C(-\beta  p- \gamma x))\geq 1$
equals $$c(d, g, \gamma):=\bigl(\gamma^2(2d-g)-\gamma\bigr){g\choose
d}.$$
\newline
\noindent $\bullet$ The number of pairs $(L, x)\in W^1_d(C)\times C$
satisfying the conditions
$$h^0(L\otimes \OO_C(-(2d-g-2) p))\geq 1 \mbox{ and }\
h^0(L\otimes \OO_C(-3 x))\geq 1 \mbox{ equals }
$$
$$e(d, g):=8\frac{g!}{(d-3)!\ (g-d+2)!}-8\frac{g!}{d!\ (g-d-1)!} .$$
\end{proposition}

We now prove more specialized results, adapted to our situation of
counting pencils with two triple points:

\begin{proposition}\label{triple1}
(1) We fix $d\geq 3$ and a general $2$-pointed curve $[C, p, q]\in
\cM_{2d-6}$. The number of pencils $l\in G^1_d(C)$ having triple
points at both $p$ and $q$ equals
$$F(d):=(2d-6)!\ \Bigl (\frac{1}{(d-3)!^2}-\frac{1}{d!\ (d-6)!}\Bigr).$$
\noindent (2) For a general curve $[C]\in \cM_{2d-4}$,  the number
of pencils $l\in G^1_d(C)$ having triple ramification at unspecified
distinct points $x,y\in C$, equals
$$N(d):=\frac{48(6d^2-28d+35)\ (2d-4)!}{d!\ (d-3)!}\mbox{ }.$$
\noindent (3) We fix a general pointed curve $[C, p]\in \cM_{2d-5,
1}$. The number of pencils $L\in W^1_d(C)$ satisfying the conditions
$$h^0(L\otimes \OO_C(-2 p))\geq 1, \ h^0(L\otimes \OO_C(-3
x))\geq 1, \ h^0(L\otimes \OO_C(-3 y))\geq 1$$ for unspecified
distinct points $x, y\in C$, is equal to
$$N_1(d):=24(12d^3-92d^2+240d-215)\frac{(2d-4)!}{d!\ (d-2)!}.$$
\end{proposition}

\begin{remark} In the formulas for $e(d, g)$ and $F(d)$ we set $1/n!:=0$ for
$n<0$.
\end{remark}
\begin{remark}As a check, for $d=3$ Proposition \ref{triple1} (2) reads
$N(3)=80$. Thus for a general curve $[C]\in \cM_2$ there are
$160=2\cdot 80$ pairs of points $(x,y)\in C\times C$, $x\neq y$,
such that $3x\equiv 3y$. This can be seen directly by considering
the map $\psi:C\times C\rightarrow \mbox{Pic}^0(C)$ given by
$\psi(x,y):=\mathcal{O}_C(3x-3y)$. Then
$\psi^*(0)=\frac{1}{2}\int_{C\times C}\psi^*(\omega \wedge
\omega)=2\cdot 3^2\cdot 3^2=162$, where $\omega$ is a differential
form representing $\theta$. To get the answer to our question we
subtract from $162$ the contribution of the diagonal
$\Delta\subseteq C\times C$. This excess intersection contribution
is equal to $2$ (cf. [Di]), so in the end we get $160=162-2$ pairs
of distinct points $(x,y)\in C\times C$ with $3x\equiv 3y$.
\end{remark}
\begin{proof} {\bf{(1)}} This is a standard exercise in limit linear series
and Schubert calculus in the spirit of \cite{EH}. We let $[C,p,q]\in
\cM_{2d-6, 2}$ degenerate to the stable $2$-pointed curve
$[C_0:=\PP^1\cup E_1\cup \ldots \cup E_{2d-6}, p_0,q_0]$, consisting
of elliptic tails $\{E_i\}_{i=1}^{2d-6}$ and a rational spine, such
that $\{p_i\}=E_i\cap \PP^1$, and the marked points $p_0, q_0$ lie
on the spine. We also assume that $p_1,\ldots,p_{2d-6},p_0,q_0\in
\mathbb P^1$ are general points, in particular $p_0, q_0\in
\PP^1-\{p_1, \ldots, p_{2d-6}\}$. Then $F(d)$ is  the number of
limit $\mathfrak g^1_d$'s on $C_0$ having triple ramification at
both $p_0$ and $q_0$ and this is the same as the number of
$\mathfrak g^1_d$'s on $\PP^1$ having cusps  at
$p_1,\ldots,p_{2d-6}$ and triple ramification at $p_0$ and $q_0$.
This  equals the intersection number of Schubert cycles
$\sigma_{(0,2)}^2 \sigma_{(0,1)}^{2d-6}$ (computed in
$H^{\mathrm{top}}(\mathbb G(1,d),\mathbb Z))$. The product can be
computed using formula (v) on page 273 in \cite{Fu} and one finds
that
$$\sigma_{(0,2)}^2\ \sigma_{(0,1)}^{2d-6}=(2d-6)!\ \Bigl(\frac{1}{(d-3)!^2}-\frac{1}{d!\ (d-6)!}\Bigr).$$

\noindent {\bf{(2)}} This is more involved. We specialize $[C]\in
\cM_{2d-4}$ to $[C_0:=\PP^1\cup E_1\cup\ldots \cup E_{2d-4}]$, where
$E_i$ are general elliptic curves, $\{p_i\}=\PP^1\cap E_i$ and
$p_1,\ldots,p_{2d-4}\in \PP^1$ are general points. Then $N(d)$ is
equal to the number of limit $\mathfrak g^1_d$'s on $C_0$ with
triple ramification at two distinct points $x,y\in C_0$. Let $l$ be
such a limit $\mathfrak g^1_d$. We can assume that both $x$ and $y$
are smooth points of $C_0$ and by the additivity of the
Brill-Noether number (see e.g. \cite{EH} pg. 365), we find that
$x,y$ must lie on the tails $E_i$. Since $[E_i, p_i]\in \cM_{1, 1}$
is general, we assume that $j(E_i)\neq 0$ (that is, none of the
$E_i$'s is the Fermat cubic). Then there can be no $l_i\in
G^1_3(E_i)$ carrying $3$ triple ramification points. There are two
cases we consider:

\noindent a) There are indices $1\leq i<j\leq 2d-4$ such that $x\in
E_i$ and $y\in E_j$. Then
$a^{l_{E_i}}(p_i)=a^{l_{E_j}}(p_j)=(d-3,d)$, hence $3x\equiv 3p_i$
on $E_i$ and $3y\equiv 3p_j$ on $E_j$. There are $8$ choices for
$x\in E_i$, $8$ choices for $y\in E_j$ and ${2d-4\choose 2}$ choices
for the tails $E_i$ and $E_j$ containing the triple points. On
$\PP^1$ we count $\mathfrak g^1_d$'s with cusps at
$\{p_1,\ldots,p_{2d-4}\}-\{p_i,p_j\}$ and triple points at $p_i$ and
$p_j$. This number is again equal to $\sigma_{(0,2)}^2\
\sigma_{(0,1)}^{2d-6}\in H^{\mathrm{top}}(\mathbb G(1, d), \mathbb
Z)$ and  we get a contribution of
\begin{equation}\label{cont1}
64{2d-4\choose 2}\sigma_{(0,2)}^2\ \sigma_{(0,1)}^{2d-6}=32(2d-4)!\ \Bigl(\frac{1}{(d-3)!^2}
-\frac{1}{d!\ (d-6)!}\Bigr)\mbox{ }.
\end{equation}
b) There is $1\leq i\leq 2d-4$ such that $x,y\in E_i$. We
distinguish between two subcases:

\noindent $b_1)$ $a^{l_{E_i}}(p_i)=(d-3,d-1)$. On $\mathbb P^1$ we
count $\mathfrak g^1_{d-1}$'s with cusps at $p_1,\ldots,p_{2d-4}$
and this number is $\sigma_{(0,1)}^{2d-4}$ (in $H^{top}(\mathbb
G(1,d-1),\mathbb Z))$. On $E_i$ we compute the number of $\mathfrak
g^1_3$'s having triple ramification at unspecified points $x,y\in
E_i-\{p_i\}$ and ordinary ramification at $p_i$. For simplicity we
set $[E_i,p_i]:=[E,p]$. If we regard $p\in E$ as the origin of $E$,
then the translation map $(x,y)\mapsto (y-x,-x)$ establishes a
bijection between the set of pairs $(x,y)\in E\times E-\Delta$,
$x\neq p\neq y\neq x$, such that there is a $\mathfrak g^1_3$ in
which $x,y,p$ appear with multiplicities $3,3$ and
 $2$ respectively, and the set of pairs $(u,v)\in E\times E-\Delta$, with $u\neq p\neq v\neq u$ such that there
is a $\mathfrak g^1_3$ in which $u,v,p$ appear with multiplicities
$3,2$ and $3$ respectively. The latter set has cardinality $16$,
hence the number of pencils $\mathfrak g^1_3$ we are counting is
$8=16/2$. All in all, we find  a contribution of
\begin{equation}\label{cont2}
8(2d-4)\ \sigma_{(0,1)}^{2d-4}=16{2d-4\choose d-1}\ .
\end{equation}

\noindent $b_2)$ $a^{l_{E_i}}(p_i)=(d-4,d)$. This time, on $\PP^1$
we look at $\mathfrak g^1_d$'s with cusps at
$\{p_1,\ldots,p_{2d-4}\}-\{p_i\}$ and a $4$-fold point at $p_1$.
Their number is $\sigma_{(0,3)}\ \sigma_{(0,1)}^{2d-5} \in
H^{\mathrm{top}}(\mathbb G(1,d),\mathbb Z))$. On $E_i$ we compute
the number of $\mathfrak g^1_4$'s for which there are distinct
points $x,y\in E_i-\{p_i\}$ such that $p_i, x, y$ appear with
multiplicities $4, 3$ and $3$ respectively. Again we set
$[E_i,p_i]:=[E,p]$ and denote by $\Sigma$ the closure in $E\times E$
of the locus
$$\{(u,v)\in E\times E-\Delta:\exists l\in G^1_4(E)\mbox{ such that }a_1^l(p)=4,\ a_1^l(u)\geq 3,\
a_1^l(v)\geq 2\}.$$ The class of the curve $\Sigma$ can be computed
easily. If $F_i$ denotes the numerical equivalence class of a fibre
of the projection $\pi_i:E\times E\rightarrow E$ for $i=1,2$, then
\begin{equation}\label{classofcorr}
\Sigma\equiv 10F_1+5F_2-2\Delta.
\end{equation}
The coefficients in this expression are determined by intersecting
$\Sigma$ with $\Delta$ and the fibres of $\pi_i$. First, one has
that $\Sigma\cap \Delta=\{(x,x)\in E\times E:x\neq p, \ 4p\equiv
4x\}$ and then $\Sigma\cap \pi_2^{-1}(p)=\{(y,p)\in E\times E:y\neq
p, \ 3p\equiv 3y\}.$ These intersections are all transversal, hence
$\Sigma \cdot \Delta=15, \Sigma \cdot F_2=8$, whereas obviously
$\Sigma \cdot F_1=3$. This proves (\ref{classofcorr}).

The number of pencils $l\subseteq |\OO_E(4p)|$ having two extra
triple points will then be equal to $1/2\ \# (\mbox{ramification
points of }\pi_2:\Sigma \rightarrow E)=\Sigma^2/2=20.$ We have
obtained in this case a contribution of
\begin{equation}\label{cont3}
20(2d-4)\ \sigma_{(0,3)}\ \sigma_{(0,1)}^{2d-5}=80{2d-4\choose d}.
\end{equation}
Adding together (\ref{cont1}),(\ref{cont2}) and (\ref{cont3}), we
obtain the stated formula for $N(d)$.

\noindent {\bf{(3)}} We relate $N_1(d)$ to $N(d)$ by specializing
the general curve from $\cM_{2d-4}$ to $[C\cup_p E]\in
\Delta_1\subset \mm_{2d-4}$, where $[C, p]\in \cM_{2d-5, 1}$ and
$[E, p]\in \mm_{1, 1}$. Under this degeneration $N(d)$ becomes the
number of admissible coverings $f:X\stackrel{d:1}\rightarrow R$
having as source a nodal curve $X$ stably equivalent to $C\cup_p E$
and as target a genus $0$ nodal curve $R$. Moreover, $f$ possesses
distinct unspecified triple ramification points $x, y\in
X_{\mathrm{reg}}$. There are a number of cases depending on the
position of $x$ and $y$.

\noindent $(3_a)$ \ $x, y\in C-\{p\}$. In this case
$\mbox{deg}(f_{C})=d$ and because of the generality of $[C, p]$,
$f_{C}$ has to be one of the finitely many $\mathfrak g^1_d$'s
having two distinct triple points and a simple ramification point at
$p\in C$. The number of such coverings is precisely $N_1(d)$. By the
compatibility condition on ramification indices at $p$, we find that
$\mbox{deg}(f_{E})=2$ and the $E$-aspect of $f$ is induced by
$|\OO_E(2p)|$. The curve $X$ is obtained from $C\cup_p E$ by
inserting $d-2$ copies of $\PP^1$ at the points in
$f_{C}^{-1}(f(p))-\{p\}$. We then map these rational curves
isomorphically to $f(E)$. This admissible cover has no automorphisms
and it should be counted with multiplicity $1$.

\noindent $(3_b)$ \ $x, y\in E-\{p\}$. The curve $[C]\in \cM_{2d-5}$
being Brill-Noether general, it carries no linear series $\mathfrak
g^1_{d-2}$, hence $\mathrm{deg}(f_{C})\geq d-1$. We distinguish two
subcases:

If $\mathrm{deg}(f_{C})=d-1$, then $f_{C}$ is one of the $a(d-1,
2d-5)$ linear series $\mathfrak g^1_{d-1}$ on $C$ having $p$ as an
ordinary ramification point. Since $C$ and $E$ meet only at $p$, we
have that $\mathrm{deg}(f_{E})=3$, and $f_E$ corresponds to a
$\mathfrak g^1_3$ on $E$ having two unspecified triple points and a
simple ramification point at $p$. There are $8$ such $\mathfrak
g^1_3$'s on $E$ (see the proof of Proposition \ref{triple1}). To
obtain a degree $d$ admissible covering, we first attach a copy
$(\PP^1)_1$ of $\PP^1$ to $E$ at the point $q\in
f_E^{-1}(f(p))-\{p\}$, then map $(\PP^1)_1$ and $C$ map to the same
component of $R$. Then we insert $d-2$ copies of $\PP^1$ at the
points lying in the same fibre of $f_{C}$ as $p$. All these rational
curves map to the same copy of $R$ as $E$. Each of these $8a(d-1,
2d-5)$ admissible coverings is counted with multiplicity $1$.

If $\mathrm{deg}(f_{C})=d$, then $f_{C}$ corresponds to one of the
$a(d, 2d-5)$ linear series $\mathfrak g^1_d$ with a $4$-fold point
at  $p$. By compatibility, $f_E$ corresponds to a $\mathfrak g^1_4$
in which $p$ and two unspecified points $x, y\in E$ appear with
multiplicities $4, 3$ and $3$ respectively. There are $20$ such
$\mathfrak g^1_4$'s on $E$, hence $20a(d, 2d-5)$ admissible
coverings.

\noindent $(3_c)$ \ $x\in E-\{p\}, y\in C-\{p\}$. In this situation
$\mathrm{deg}(f_{C})=d$ and $f_{C}$ corresponds to one of the $e(d,
2d-5)$ coverings $\mathfrak g^1_d$ on $C$ having a triple point at
$p$ and another unspecified triple point at $y\in C$. Then
$\mbox{deg}(f_{E})=3$ and $3x\equiv 3p$, that is, there are $8$
choices of the $E$-aspect of $f$. We obtain $X$ by attaching to $C$
copies of $\PP^1$ at the $d-3$ points in $f_C^{-1}(f(p))-\{p\}$, and
mapping these curves isomorphically onto $f(C)$.

By degeneration to $[C\cup_p E]$, we have found the relation for
$[C, p]\in \cM_{2d-5, 1}$:
$$N(d)=N_1(d)+ 20a(d, 2d-5)+ 8a(d-1, 2d-5)+ 8e(d, 2d-5).$$ This
immediately leads to the claimed expression for $N_1(d)$.
\end{proof}

\section{The class of the divisor $\overline{\mathfrak{TR}}_d$}

The strategy to compute the class $[\overline{\mathfrak{TR}}_d]$ is
similar to the one employed by Eisenbud and Harris in \cite{EH} to
determine the class of the Brill-Noether divisors $[\mm_{g, d}^r]$
of curves with a $\mathfrak g^r_d$ in the case $\rho(g, r, d)=-1$:
We determine the restrictions of $\tr_d$ to $\mm_{0, g}$ and
$\mm_{2, 1}$ via obvious flag maps. However, because in the
definition of $\tr_d$ we allow $2$ degrees of freedom for the triple
ramification points, the calculations are much more intricate (and
interesting) than in the case of Brill-Noether divisors.

\begin{proposition}\label{flag}
Consider the flag map $j:\mm_{0, g}\rightarrow \mm_{g}$ obtained by
attaching $g$ general elliptic tails at the $g$ marked points. Then
$j^*(\overline{\mathfrak{TR}}_d)=0$. If we have a linear relation
$$\tr_d\equiv a\ \lambda-\sum_{i=0}^{d-2} b_i\ \delta_i\in \mathrm{Pic}(\mm_g), \mbox{  then }
b_i=\frac{i(g-i)}{g-1} b_1, \ \mbox{ for }\ 1\leq i\leq d-2.$$
\end{proposition}
\begin{proof} The second part of the statement is a
consequence of the first: For an effective divisor $D\equiv a
\lambda-\sum_{i=0}^{d-2} b_i \delta_i\in \mathrm{Pic}(\mm_g)$
satisfying the condition $j^*(D)=\emptyset$, we have the  relations
among its coefficients: $b_i=\frac{i(g-i)}{g-1} b_1$ for $i\geq 1$
(cf. \cite{EH} Theorem 3.1).

Suppose that $[X:=R\cup_{x_1} E_1\cup \ldots \cup_{x_g} E_g]\in
j(\mm_{0, g})$ is a flag curve corresponding to a $g$-stable
rational curve $[R, x_1, \ldots, x_g]$. The elliptic tails
$\{E_i\}_{i=1}^{g}$ are general and we may assume that all the
$j$-invariants are different from $0$. In particular, none of the
$[E_i, x_i]$'s carries a $\mathfrak g^1_3$ with triple ramification
points at $x_i$ and at two unspecified points $x, y\in E_i-\{x_i\}$.
Assuming that $[X]\in \tr_d$, there exists $l\in
\overline{G}^1_d(X)$ a limit $\mathfrak g^1_d$, together with
distinct ramification points $x\neq y\in X$, such that $a_1^l(x)\geq
3$ and $a^l_1(y)\geq 3$. By blowing-up if necessary the nodes $x_i$
(that is, by inserting chains of $\PP^1$'s at the points $x_i$), we
may assume that both $x, y$ are smooth points of $X$.

We make use of the following facts: On $R$ we have that the
inequality
$$\rho(l_R, x_1, \ldots, x_g, z_1, \ldots, z_t)\geq 0,$$ for any
choice of distinct points $z_1, \ldots, z_t\in R-\{x_1, \ldots,
x_g\}$. On the elliptic tails, we have that $\rho(l_{E_i}, x_i,
z)\geq -1$, for any point $z\in E_i-\{x_i\}$, with equality only if
$z-x_i\in \mathrm{Pic}^0(E_i)$ is a torsion class. Using these
remarks as well as and the additivity of the Brill-Noether number of
$l$, since $\rho(l, x, y)=-3$ it follows that there must exist an
index $1\leq i\leq g$ such that $x, y\in E_i-\{x_i\}$, and
$\rho(l_{E_i}, x_i, x, y)=-3$. This implies that
$a^{l_{E_i}}(x_i)=(d-3, d)$ and that $l_{E_i}(-(d-3) x_i)\in
G^1_3(E_i)$ has triple ramification points at distinct points $x_i,
x$ and $y$. This can happen only if $E_i$ is isomorphic to the
Fermat cubic, a contradiction.
\end{proof}

The next result highlights the difference between $\tr$ and all the
other geometric divisors in the literature, cf. \cite{HM},
\cite{EH}, \cite{H}, \cite{Fa}, \cite{Fa2}: $\tr$ is the first
example of a geometric divisor on $\mm_g$ not pulled-back from the
space $\mm_g^{\mathrm{ps}}$ of pseudo-stable curves.

\begin{proposition}\label{elltails}
If $\tr_d\equiv a\ \lambda-\sum_{i=0}^{d-2} b_i \ \delta_i\in
\mathrm{Pic}(\mm_g)$, then $a-12b_0+b_1=4a(d, 2d-4)$.
\end{proposition}
\begin{proof}
We use a standard test curve in $\mm_g$ obtained by attaching to the
marked point of a general pointed curve $[C, q]\in \cM_{2d-4, 1}$ a
pencil of plane cubics. If $R\subset \mm_g$ is the family induced by
this pencils, then clearly $R\cdot \lambda=1, R\cdot \delta_0=12,
R\cdot \delta_1=-1$ and $R\cdot \delta_j=0$ for $j\geq 2$.

Set-theoretically,  $R\cap \tr_d$ consists of the points
corresponding to the elliptic curves $[E, q]$ in the pencil, for
which there exists $l\in G^1_3(E)$ as well as two distinct points
$x, y\in E-\{q\}$ with $a_1^l(q)=a_1^l(x)=a_1^l(y)=3$ (It is a
standard limit linear series argument to show that the triple points
of the limit $\mathfrak g^1_d$ must specialize to the elliptic
tail). Then $E$ must be isomorphic to the Fermat cubic, (thus
$j(E)=0$, and this curve appears $12$ times in the pencil. The
pencil $l\in G^1_3(E)$ is of course uniquely determined. Since $
\mathrm{Aut}(E, q)=\mathbb Z/6\mathbb Z$ while a generic element
from $\mm_{1, 1}$ has automorphism group $\mathbb Z/2\mathbb Z$,
each point of intersection will contribute $4=24/6$ times in the
intersection $R\cap \tr_d$. On the side of the genus $2d-4$
component, we count pencils $L\in W^1_d(C)$ with $a_1^L(q)\geq 3$.
Using Proposition \ref{harris}  their number is finite and equal to
$a(d, 2d-4)$, hence $R\cdot \tr_d=4a(d, 2d-4)$.
\end{proof}

Next we describe the restriction of $\tr_d$ under the map
$\chi:\mm_{2, 1}\rightarrow \mm_{2d-3}$ obtained by attaching a
fixed tail $B$ of genus $2d-5$ to each pointed curve $[C, p]\in
\cM_{2, 1}$. It is revealing to compare Theorem \ref{genus2} to
Propositions 4.1 and 5.5 in \cite{EH}: When $\rho(g, r, d)=-1$, the
pull-back of the Brill-Noether divisor $\chi^*(\mm_{g, d}^r)$ is
irreducible and supported on $\ww$. By contrast, $\tr_d$ displays a
much richer geometry.

\noindent \emph{Proof of Theorem \ref{genus2}.} We fix  a general
pointed curve $[B, p]\in \cM_{2d-5, 1}$. For each $[C, p]\in \cM_{2,
1}$, we study degree $d$ admissible coverings $[f:X\rightarrow R,
q_1, q_2; p_1, \ldots, p_{6d-12}]\in \hh_d$ with source curve $X$
stably equivalent to $C\cup_p B$, and target $R$ a nodal curve of
genus $0$. Moreover, $f$ is assumed to have distinct points of
triple ramification $x, y\in X_{\mathrm{reg}}$, where $f(x)=q_1$ and
$f(y)=q_2$. It is easy to check that both $x$ and $y$ must lie
either on $C$ or on $B$ (and not on rational components of $X$ we
may insert). Depending on their position we distinguish four cases:

\noindent (i) $x, y\in B$. A parameter count shows that
$\mbox{deg}(f_{B})=d$ and $p\in B$ must be a simple ramification
point for $f_{B}$. By compatibility of ramification sequences at
$p$, then $f_{C}$ must also be simply ramified at $p$, that is,
$p\in C$ is a Weierstrass point and $f_C$  is induced by
$|\OO_C(2p)|$.  There is a canonical way of completing $\{f_{C},
f_{B}\}$ to an element  in $\mathfrak{H}_d$, by attaching rational
curves to $B$ at the points in $f_{B}^{-1}(f(p))-\{p\}$. For a fixed
$[C, p]\in \ww$, the Hurwitz scheme is smooth at each of the points
$t\in \hh_d$ corresponding to an admissible coverings $\{f_{C},
f_{B}\}$ of the type described above. Since $t$ has no automorphisms
permuting some of the branch points, it follows that
$\mathfrak{H}_d=\hh_d/\mathfrak S_2\times \mathfrak S_{6d-12}$ is
also smooth at each of the $N_1(d)$ points in the fibre
$\sigma^{-1}([C\cup_p B])$. This implies that $N_1(d)\cdot \ww$
appears as an irreducible component in the pull-back divisor
$\chi^*(\tr_d)$.

\noindent (ii) $x, y\in C, \mbox{ deg}(f_{B})=d$. Clearly
$\mbox{deg}(f_{C})\geq 4$ and the $B$-aspect of the covering must
have a $4$-fold point at $p$. There are $a(d, 2d-5)$ choices for
$f_{B}$, whereas $f_{C}$ corresponds to a linear series $l_C\in
G^1_4(C)$ with $a_1^{l_C}(p)=4$ and which has two other points of
triple ramification. To obtain the domain of an admissible covering,
we attach to $B$ rational curves at the $(d-4)$ points in
$f_{B}^{-1}(f(p))-\{p\}$. We map these curves isomorphically onto
$f_C(C)$. The divisor $a(d, 2d-5)\cdot \dd_3$ is an irreducible
component of $\chi^*(\tr_d)$.

\noindent (iii) $x, y\in C,\  \mbox{deg}(f_{B})=d-1$. In this case
the $B$-aspect corresponds to one of the $a(d-1, 2d-5)$ linear
series $l_B\in G^1_{d-1}(B)$ with simple ramification at $p$, while
$f_C$ is a degree $3$ covering having two unspecified points of
triple ramification and simple ramification at $p\in C$. To obtain a
point in $\mathfrak{H}_d$, we attach a rational curve $T'$ to $C$ at
the remaining point in $f_{C}^{-1}(f(p)-\{p\}$. We then map $T'$
isomorphically onto $f_B(B)$. Next, we attach $d-3$ rational curves
to $B$ at the points $f_{B}^{-1}(f(p))-\{p\}$, which we map
isomorphically onto $f_C(C)$. Each resulting admissible covering has
no automorphisms and is a smooth point of $\mathfrak{H}_d$. Thus
$a(d-1, 2d-5)\cdot \overline{\mathcal{D}}_2$ is a component of
$\chi^*(\tr_d)$.

\noindent (iv) $x\in C, y\in B$. After a moment of reflection we
conclude that $\mbox{deg}(f_{B})=d$, that is,  $f_{B}$ corresponds
to one of the $e(d, 2d-5)$ coverings $l_B\in G^1_d(B)$ with
$a_1^{l_B}(p)=3$  and $a_1^{l_B}(y)=3$ at some unspecified point
$y\in B-\{p\}$. The $C$-aspect of $f$ is determined by the choice of
a point $x\in C-\{p\}$ such that $3x\equiv 3p$. Hence $e(d,
2d-5)\cdot \overline{\mathcal{D}}_1$ is the final irreducible
component of $\chi^*(\tr_d)$. \hfill $\Box$

As a consequence of Proposition \ref{flag} and Theorem \ref{genus2}
we are in a position to determine all the $\delta_i$-coefficients
($i\geq 1$) in the expansion of $\tr_d$ in the basis of
$\mbox{Pic}(\mm_g)$:

\begin{theorem}\label{higherdeltas}
If $\tr_d \equiv a\ \lambda-\sum_{i=0}^{d-2} b_i\ \delta_i \in
\mathrm{Pic}(\mm_g)$, then we have that
$$b_i=\frac{(2d-6)!}{2\ d! (d-3)!} \ i(2d-3-i)(36d^3-156d^2+180d-5),
\ \mbox{ for all }1\leq i\leq d-2.$$
\end{theorem}
\begin{proof}
We use the obvious relations $\chi^*(\delta_2)=-\psi$,
$\chi^*(\lambda)=\lambda, \ \chi^*(\delta_0)=\delta_0, \
\chi^*(\delta_1)=\delta_1$. If for a class $E\in
\mathrm{Pic}(\mm_{2, 1})$ we denote by $(E)_{\psi}$ the coefficient
of $\psi$ in its expansion in the basis $\{\psi, \lambda,
\delta_0\}$ of $\mathrm{Pic}(\mm_{2, 1})$ (see also the next section
for details on the divisor theory of $\mm_{2, 1}$), then, using
Proposition \ref{elltails}, we can write the following relation:
$$b_2=\frac{2(g-2)}{g-1} b_1=N_1(d)(\ww)_{\psi} +e(d, 2d-5)
(\overline{\mathcal{D}}_1)_{\psi}+a(d-1, 2d-5)
(\overline{\mathcal{D}}_2)_{\psi}+a(d, 2d-5)(
\overline{\mathcal{D}}_3)_{\psi}.$$ We determine the coefficients
$(\dd_i)_{\psi}$ for $1\leq i\leq 3$ by intersecting each of these
divisors with a general fibral curve $F:=\{[C, p]\}_{p\in C}\subset
\mm_{2, 1}$ of the projection $\pi:\mm_{2, 1}\rightarrow \mm_2$.
(Note that $(\ww)_{\psi}=3$).

It is useful to recall that if $[C, q]\in \cM_{2, 1}$ is a fixed
general pointed curve and $a\geq b\geq 0$ are integers, then the
number of pairs $(p, x)\in C\times C, p\neq x$ satisfying a linear
equivalence relation $a\cdot x\equiv b\cdot p+(a-b)\cdot q$ in
$\mbox{Pic}^a(C)$, equals
\begin{equation}\label{r}
r(a, b):=2(a^2 b^2-1).
\end{equation}
We start with $\dd_1$ and note that $F\cdot \dd_1$ is the number of
pairs $(x, p)\in C\times C$ with $x\neq p$, such that $3x\equiv 3p$,
which is equal to $r(3, 3)=160$ and then $(\dd_1)_{\psi}=r(3,
3)/(2g-2)=80$. To compute $F\cdot \dd_2$ we note that there are
$80=r(3, 3)/2$ pencils $L\in W^1_3(C)$ with two distinct triple
ramification points. From the Hurwitz-Zeuthen formula, each such
pencil has $4$ more simple ramification points, thus
$(\dd_2)_{\psi}=4\times 80/(2g-2)=160$. Finally, $F\cdot
\dd_3=n_0/2$, where by $n_0$ we denote the number of pencils $l\in
W^1_4(C)$ having one unspecified point of total ramification and two
further points of triple ramification, that is there exist mutually
distinct points $x, y, p\in C$ with $a_1^l(p)=4$ and
$a_1^l(x)=a_1^l(y)=3$.

We compute $n_0$ by letting $C$ specialize to a curve of compact
type $[C_0:=C_1\cup_q C_2]$, where $[C_1, q], [C_2, q]\in \cM_{1,
1}$. Then $n_0$ is the number of admissible coverings
$f:X\stackrel{4:1} \rightarrow R$, where $R$ is of genus $0$ and $X$
is stably equivalent to $C_0$ and has a $4$-fold ramification point
$p\in X_{\mathrm{reg}}$ and triple ramification  points $x, y\in
X_{\mathrm{reg}}$. We distinguish three cases:

\noindent (i) $x, y\in C_2$ and $p\in C_1$ (Or $x, y\in C_1$ and
$p\in C_2$). In this case
$\mbox{deg}(f_{C_1})=\mbox{deg}(f_{C_2})=4$ and we have the linear
equivalence $4p\equiv 4q$ on $C_1$. This yields $15$ choices for
$p\neq q$. On $C_2$ we count $\mathfrak g^1_4$'s with total
ramification at $q$, and two unspecified triple points. This number
is equal to $20$ (see the proof of Proposition \ref{triple1}).
Reversing the role of $C_1$ and $C_2$ we double the number of
coverings and we find $600=2\cdot 15\cdot 20$ admissible $\mathfrak
g^1_4$'s.

\noindent (ii) $x, p\in C_2$ and $y\in C_1$ (Or $x, p\in C_1$ and
$y\in C_2$). In this situation $\mbox{deg}(f_{C_1})=3$ and
$\mbox{deg}(f_{C_2})=4$ and on $C_1$ we have the linear equivalence
$3y\equiv 3q$, which gives $8$ choices for $y$. On $C_2$ we count
$l_{C_2}\in G^1_4(C_2)$ in which two unspecified points $p, x\in
C_2$ appear with multiplicities $4$ and $3$ respectively, while
$a_1^{l_{C_2}}(q)=3$. By translation, this is the same as the number
of pairs of distinct points $(u, v)\in C_2-\{q\}\times C_2-\{q\}$
such that there exists $l_2\in G^1_4(C_2)$ with $a_1^{l_2}(q)=4,
a_1^{l_2}(x)=a_1^{l_2}(y)=3$. This number equals $40$ (again, see
the proof of Proposition \ref{triple1}). By reversing the role of
$C_1$ and $C_2$ the total number of coverings in case (ii) is $640=2
\cdot 8\cdot 40$.

\noindent (iii) $x, y, p\in C_1$ (or $x, y, p\in C_2$). A quick
parameter count shows that $\mbox{deg}(f_{C_2})=2$ and
$\mbox{mult}_q(f_{C_2})=\mbox{mult}_q(f_{C_1})=2$. Hence $f_{C_2}$
is induced by $|\OO_{C_2}(2q)|$. On $C_1$ we count $\mathfrak
g^1_4$'s in which the points $p, x, y, q$ appear with multiplicities
$4, 3, 3$ and $2$ respectively. The translation on $C_2$ from $p$ to
$q$ shows that we are yet again in the situation of Proposition
\ref{triple1} and this last number is $20$. We interchange $C_1$ and
$C_2$ and we find $40$ admissible $\mathfrak g^1_4$'s on $C_1\cup
C_2$ with all the non-ordinary ramification concentrated on a single
component.

By  adding (i), (ii) and (iii) together, we obtain
$n_0=600+640+40=1280$. This determines
$(\dd_3)_{\psi}=n_0/(2g-2)=640$ and completes the proof.
\end{proof}

\section{The divisor theory of $\overline{\mathcal{M}}_{2,1}$}
The remaining part of the calculation of $[\tr_d]$ has been reduced
to the problem of determining the divisor classes
$[\overline{\mathcal{D}}_i]$ $(i=1, 2, 3)$ on $\mm_{2, 1}$. We
recall some things about divisor theory on this space (see also
\cite{EH}). There are two boundary divisor classes:

\noindent $\bullet$ $\delta_0$, whose generic point is an
irreducible $1$-pointed nodal curve of genus $2$.
\newline
\noindent $\bullet$ $\delta_1$, with generic point being a
transversal union of two elliptic curves with the marked point lying
on one of the components.

If $\pi:\mm_{2, 1}\rightarrow \mm_2$ is the universal curve then
$\psi:=c_1(\omega_{\pi})\in \mathrm{Pic}(\mm_{2, 1})$ denotes the
tautological class and $\lambda=\pi^*(\lambda) \in
\mathrm{Pic}(\mm_{2, 1})$ is the Hodge class. Unlike the case $g\geq
3$, $\lambda$  is a boundary class on $\mm_2$, and we have Mumford's
genus $2$ relation:
$$\lambda=\frac{1}{10}\delta_0+\frac{1}{5}\delta_1.$$
The classes $\psi, \lambda$ and $\delta_1$ form a basis of
$\mbox{Pic}(\mm_{2, 1})\otimes \mathbb Q$. The class of the
Weierstrass divisor has been computed in \cite{EH} Theorem 2:
\begin{equation}\label{weierstrass}
\ww\equiv 3 \psi-\lambda-\delta_1.
\end{equation}

We start by determining the class of $\overline{\mathcal{D}}_1$ of
$3$-torsion points:
\begin{proposition}\label{d1}
The class of the closure in $\mm_{2, 1}$ of the effective divisor
$$\mathcal{D}_1=\{[C,p]\in \mathcal{M}_{2,1}:\exists x\in C-\{p\}\ \mbox{ such that }\ 3x\equiv 3p\}$$
is given by $[\overline{\mathcal{D}}_1]=80 \psi+10 \delta_0-120
\lambda\in \mathrm{Pic}(\mm_{2, 1})$.
\end{proposition}
\begin{proof} We introduce the map $\chi:\overline{\mathcal{M}}_{2,1}\rightarrow
\overline{\mathcal{M}}_4$ given by $\chi([C,p]):=[B\cup_p C]$, where
$[B,p]$ is a general $1$-pointed curve of genus $2$. On
$\overline{\mathcal{M}}_4$ we have the divisor of curves with an
exceptional Weierstrass point $\mathfrak{Di}:=\{[C]\in
\mathcal{M}_4:\exists x\in C\mbox{ such that }h^0(C,3x)\geq 2\}.$
Its class has been computed by Diaz \cite{Di}:
$\overline{\mathfrak{Di}}\equiv
264\lambda-30\delta_0-96\delta_1-128\delta_2 \in
\mathrm{Pic}(\mm_4)$.

We claim that $\chi^*(\overline{\mathfrak{Di}})=\dd_1+16\cdot
\overline{\mathcal{W}}$. Indeed, let $[C,p]\in \mathcal{M}_{2,1}$ be
such that $\chi([C,p])\in \overline{\mathfrak{Di}}$. Then there is a
limit $\mathfrak g^1_3$ on $X:=B\cup_p C$, say $l=\{l_B, l_C\}$,
which has a point of total ramification at some $x\in
X_{\mathrm{reg}}$. There are two possibilities:

\noindent  (i) If $x\in C$, then $a^{l_{B}}(p)=(0,3)$, hence
$l_{B}=|\OO_B(3p)|$, while on  $C$ we have the linear equivalence
$3p\equiv 3x$, that is, $[C,p]\in \dd_1$.

\noindent (ii) If $x\in B$, then $a^{l_C}(p)=(1,3)$, that is, $p\in
B$ is a Weierstrass point and moreover $l_C=p+|\OO_C(2p)|$. On $B$
we have that $a^{l_{B}}(p)=(0,2)$ and $a^{l_{B}}(x)=(0,3)$, that is,
$3x\equiv 2p+y$ for some $y\in B-\{p, y\}$. There are $r(3, 1)=16$
such pairs $(x, y)$.

Thus we have proved that
$\chi^*(\overline{\mathfrak{Di}})=\dd_1+16\cdot
\overline{\mathcal{W}}$ (We would have obtained the same conclusion
using admissible coverings instead of limit $\mathfrak g^1_3$'s). We
find the formula for $[\dd_1]$ if we remember that
$\chi^*(\delta_0)=\delta_0,\mbox{ }\chi^*(\delta_1)=\delta_1,\mbox{
}\chi^*(\delta_2)=-\psi$ and $\chi^*(\lambda)=\lambda.$
\end{proof}

\subsection{The divisor $\tr_3$  and the class of $\dd_2$}
We compute the class of the divisor $\dd_2$ on $\mm_{2, 1}$ by
determining directly the class of $\tr_3$ in genus $3$ (In this case
$\dd_3=\emptyset$). Much of the set-up we develop here is valid for
arbitrary $d\geq 3$ and will be used in the next section when we
compute the class $[\tr_4]$ on $\mm_5$. We fix a general $[C, p]\in
\cM_{2d-4, 1}$  and introduce the following enumerative invariant:
$$N_2(d):=\#\{l\in G^1_d(C):\exists x\neq y\in C-\{p\} \mbox{ such
that } \ l(-3 x)\neq \emptyset \  \mbox{ and } \ l(-p-2 y)\neq
\emptyset\}.$$ For instance, $N_2(3)$ is the number of pairs $(x,
y)\in C \times C$, $x\neq p\neq y$ such that $3x\equiv p+2y$, hence
$N_2(3)=r(3, 2)=70$ (cf. formula (\ref{r})).

\noindent For each $d\geq 4$ we fix a general pointed curve $[B,
q]\in \cM_{2d-5, 1}$ and define the invariant:
$$N_3(d):=\#\{l\in G^1_d(B):\exists x\neq y\in B-\{q\} \mbox{ such
that }\ l(-3x)\neq \emptyset \mbox{ and } \ l(-2q-2y)\neq
\emptyset\}.$$

\begin{theorem}\label{genus 3}
The closure of the divisor $\mathfrak{TR}_3:=\{[C]\in \cM_3: \exists
x\neq p\in C \mbox{ with } 3x\equiv 3x\}$ is linearly equivalent to
the class
$$\tr_3 \equiv 2912 \lambda-311 \delta_0- 824 \delta_1\in \mathrm{Pic}(\mm_3).$$ It follows
that $\dd_2\equiv -200 \lambda+160 \psi+17  \delta_0\in
\mathrm{Pic}(\mm_{2, 1})$.
\end{theorem}
\begin{proof}
For most of this proof we assume $d\geq 3$ and we specialize to the
case of $\mm_3$ only at the very end. We write $\tr_d \equiv a\
\lambda-b_0\ \delta_0- \cdots -b_{d-2}\ \delta_{d-2}\in
\mbox{Pic}(\mm_g)$ and we have already determined $b_1, \ldots,
b_{d-2}$ (cf. Theorem \ref{higherdeltas}) while we know that
$a-12b_0+b_1=4a(d, 2d-4)$ (cf. Proposition \ref{elltails}). We need
one more relation involving $a, b_0$ and $b_1$, which we obtain by
intersecting $\tr_d$ with the test curve
$$C^0:=\bigl\{\frac{C}{q\sim p}\bigr\}_{p\in C}\subset \Delta_0\subset \mm_g$$
obtained from a  general curve $[C, q]\in \cM_{2d-4, 1}$. The number
$C^0\cdot \tr_d$  counts (with appropriate multiplicities)
admissible coverings $$t:=[f:X\stackrel{d:1}\rightarrow R, \ q_1,
q_2: p_1, \ldots, p_{6d-12}] \ \mathrm{mod} \mbox{ }\mathfrak
S_2\times \mathfrak S_{6d-12}\in \mathfrak{H}_d,$$ where the source
$X$ is stably equivalent to the curve $C\cup_{\{p, q\}} T$ $(q\in
 C)$ obtained by "blowing-up" $\frac{C}{q\sim p}$ at the node and inserting
 a rational curve $T$. These covers should possess two points of triple
ramification  $x, y\in X_{\mathrm{reg}}$ such that $f(x)=q_1,
f(y)=q_2$. Suppose $t\in C^0\cdot \tr$ and again we distinguish a
number of possibilities:

\noindent (i) $x, y\in C$. Then $\mbox{deg}(f_{C})=d$ and $f_{C}$
corresponds to one of the $N(d)$ linear series $l\in G^1_d(C)$ with
two points of triple ramification. The point $q\in C$ is such that
$l(-p-q)\neq \emptyset$, which, after having fixed $l$, gives $d-1$
choices. Clearly $\mbox{mult}_{q}(f_{C}) =\mbox{mult}_q(f_{T})=1$.
This implies that $\mbox{deg}(f_{T})=2$ and $f_{T}$ is given by
$|\OO_T(p+q)|$. To obtain out of $\{f_C, f_B\}$  a point $t\in
\hh_d$, we attach rational curves to $C$ at the points in
$f_{C}^{-1}(f(p))-\{p, q\}$ and map these isomorphically onto the
component $f_T(T)$ of $R$. Each such cover has an automorphism
$\phi:X\rightarrow X$ of order $2$ such that
$\phi_{C}=\mathrm{id}_C$, $\phi_{T'}=\mathrm{id}_{T'}$, for every
rational component $T'\neq T$ of $X$, but $\phi_{T}$ interchanges
the $2$ branch points of $T$. Even though $t\in \hh_d$ is a smooth
point (because there is no automorphism of $X$ preserving \emph{all}
the ramification points of $f$), if $\tau\in \mathfrak{S}_{6d-12}$
is the involution exchanging the marked points lying on $f_T(T)$,
then $\tau\cdot t=t$. Therefore $\hh_d/\mathfrak{S}_2\rightarrow
\mm_g$ is simply ramified at $t$. In a general deformation
$[\mathcal{X}\rightarrow \mathcal{R}]$ of $[f: X\rightarrow R]$ in
$\hh_d$ we blow-down  $T$ and obtain a rational double point, hence
the image of $\mathcal{R}$ in $\mm_g$ meets $\Delta_0$ with
multiplicity $2$.  Since $\hh_d/\mathfrak{S}_2\rightarrow \mm_g$ is
ramified anyway, it follows that each of the $(d-1)N(d)$ admissible
coverings found at this step is to be counted with multiplicity $1$.

\noindent (ii) $x\in C, y\in T$. Since $C$ has only finitely many
$\mathfrak g^1_{d-1}$'s, all simply ramified and having no
ramification in the fibre over $q$, we must have that
$\mbox{deg}(f_{C})=d$ and $\mbox{deg}(f_{T})=3$. Moreover, $C$ and
$T$ map via $f$ onto the two components of the target $R$ in such a
way that $f_{C}(p)=f_{C}(q)=f_{T}(p)=f_{C}(q)$. In particular, both
$f_{C}$ and $f_{T}$ are simply ramified at either $p$ or $q$. If
$f_{C}$ is ramified at $q\in C$, then $f_C$ is induced by one of the
$e(d, 2d-4)$ linear series $l\in G^1_d(C)$ with one unassigned point
of triple ramification and one assigned point of simple
ramification. Having fixed $l$, there are $d-2$ choices for $p\in C$
such that $l(-2q-p)\neq \emptyset$. On $T$ there is a unique
$\mathfrak g^1_3$ corresponding to a map $f_{T}:T\rightarrow \PP^1$
such that $f_T^*(0)=2 q+p$ and $f_T^*(\infty)=3 y$, for some $y\in
T-\{q, p\}$. Finally, we attach $d-3$ rational curves to $C$ at the
points in $f_C^{-1}(f(q))-\{p, q\}$ and we map these components
isomorphically onto $f_T(T)$.

The other possibility is that $f_{C}$ is unramified at $q$ and
ramified at $p$. The number of such $\mathfrak g^1_d$'s is $N_2(d)$.
On the side of $T$, there is a unique way of choosing
$f_{T}:T\stackrel{3:1}\rightarrow \PP^1$ such that $f_T^*(0)=q+2 p$
and $f_T^*(\infty)=3 y$. Because the map
$\sigma:\mathfrak{H}_d\rightarrow \mm_g$ blows-down the component
$T$, if $[\mathcal{X}\rightarrow \mathcal{R}]$ is a general
deformation of $[f:X\rightarrow R]$ then $\sigma(\mathcal{R})$ meets
$\Delta_0$ with multiplicity $3$ (see also \cite{Di}, pg. 47-52).
Thus $\tr_d \cdot \Delta_0$ has multiplicity $3$ at the point
$[C/p\sim q]$. The admissible coverings constructed at this step
have no automorphisms, hence they each must be counted with
multiplicity $3$. This yields a total contribution of $3(d-2)e(d,
2d-4)+3N_2(d)$.

\noindent (iii) $x, y\in T-\{p, q\}$. Here there are two subcases.
First, we assume that $\mbox{deg}(f_C)=d-1$, that is, $f_C$ is
induced by one of the $\frac{(2d-4)!}{(d-1)! (d-2)!}$ linear series
$l\in G^1_{d-1}(C)$. For each such $l$, there are $d-2$
possibilities for $p$ such that $l(-q-p)\neq \emptyset$. Clearly
$\mbox{deg}(f_{T})=3$ and the admissible covering $f$ is constructed
as follows: Choose $f_{T}:T\rightarrow \PP^1$ such that $f_T^*(0)=3
x$, $f_T^*(\infty)=3 y$ and $f_T^*(1)=p+q+q'$. We map $C$ to the
component of $R$ other than $f_T(T)$ by using $l\in G^1_{d-1}(C)$
and $f_C(p)=f_T(p)$ and $f_C(q)=f_T(q)$. We attach to $T$ a rational
curve $T'$ at the point $q'$ and map $T'$ isomorphically onto
$f(C)$. Finally we attach $d-3$ rational curves to $C$ at the points
in $f_C^{-1}(f(q))-\{q, p\}$. Each of these ${2d-4\choose d-1}$
elements of $\mathfrak h_d$ is counted with multiplicity $2$.

We finally deal with the case $\mbox{deg}(f_C)=d$. Since a
$\mathfrak g^1_3$ on $\PP^1$ with two points of total ramification
must be unramified everywhere else, it follows that
$\mbox{deg}(f_T)\geq 4$. The generality assumption on $[C, q]$
implies that $\mbox{deg}(f_T)=4$. The $C$-aspect of $f$ is induced
by $l\in G^1_d(C)$ for which there are integers $\beta, \gamma \geq
1$ with $\beta+\gamma=4$ and a point $p\in C$ such that $l(-\beta
p-\gamma q)\neq \emptyset$.  Proposition \ref{harris} gives the
number $c(d, 2d-4, \gamma)$ of such $l\in G^1_d(C)$. On the side of
$T$, we choose $f_T:T \stackrel{4:1}\rightarrow \PP^1$ such that
$f_T^*(0)=3 x,\ f_T^*(\infty)=3 y$ and $f_T^*(1)=\beta p+\gamma q$.
When $\gamma\in \{1, 3\}$, up to isomorphism there is a unique such
$f_T$ having $3$ triple ramification points. By direct computation
we have the formula:
$$f_T:T \rightarrow \PP^1, \mbox{   }\ f_T(t):=\frac{2t^3(t-2)}{2t-1},$$
which has the properties that
$f_T^{(i)}(0)=f_T^{(i)}(\infty)=f_T^{(i)}(1)=0$, for $i=1, 2$. When
$\gamma=2$, there are two $\mathfrak g^1_4$'s with $2$ points of
triple ramification and $2$ points of simple ramification lying in
the same fibre. It is important to point out that $f_T$ (and hence
the admissible covering $f$ as well), has an automorphism of order
$2$ which preserves the points of attachment $p, q\in T$ but
interchanges $x$ and $y$ (In coordinates, if $x=0, y=\infty \in T$,
check that $f_T(1/t)=1/f_T(t)$). This implies that $\hh_d\rightarrow
\mm_d$ is (simply) ramified at $[X\rightarrow R]$. Furthermore, a
calculation similar to \cite{Di} pg. 47-50, shows that the image in
$\mm_g$ of a generic deformation in $\hh_d$ of $[X\rightarrow T]$
meets the divisor $\Delta_0$ with multiplicity $4=\beta+\gamma$. It
follows that $\tr_d \cdot \Delta_0$ has multiplicity $4/2=2$ in a
neighbourhood of $[C/p\sim q]$, that is, each covering found at this
step gets counted with multiplicity $2$ in the product $C^0\cdot
\tr$. Coverings of this type give a contribution of
$$2c(d, 2d-4, 1)+2c(d, 2d-4, 3)+4c(d, 2d-4, 2)=128{2d-4\choose d}.$$
Thus we can write the following equation:
\begin{equation}\label{ycontr}(2g-2)b_0-b_1=C^0\cdot
\tr_d =\end{equation} $$=(d-1)N(d)+3N_2(d)+3(d-2)e(d,
2d-4)+128{2d-4\choose d}+2{2d-4\choose d-1}. $$ For $d=3$, when
$N_2(d)=70$, all terms in (\ref{ycontr}) are known and this finishes
the proof.
\end{proof}

\section{The divisor $\tr_5$ and the class of $\dd_3$}

In this section we finish the computation of $[\tr_d]$ (and
implicitly compute $[\dd_3]\in \mathrm{Pic}(\mm_{2, 1})$ and
determine $N_2(d)$ for all $d\geq 3$ as well). According to
(\ref{ycontr}) it suffices to compute $N_2(4)$ to determine
$[\tr_4]\in \mbox{Pic}(\mm_5)$. Then applying Theorem \ref{genus2}
we obtain $[\dd_3]$ which will finish the calculation of $[\tr_d]$
for
 $g=2d-3$. We summarize some of the enumerative results needed in
this section:

\begin{proposition}\label{elliptic}
We fix a general $2$-pointed elliptic curve $[E, p, q]\in \cM_{1,
2}$.

\noindent (a) There are $11$ pencils $l\in G^1_3(E)$ such that there
exist distinct points $x, y\in E-\{p, q\}$ with $a^l_1(x)=3,\
a^l_1(q)=2$ and $l(-p-2 y)\neq \emptyset$.

\noindent (b) There are $38$ pencils $l\in G^1_4(E)$ such that there
exist distinct points $x, y\in E-\{p, q\}$ with $a^l_1(p)=4, \
a^l_1(x)=3$  and $l(-q-2 y)\neq \emptyset$.
\end{proposition}
\begin{proof}
{\bf{(a)}} We denote by $\cU$ the closure in $E\times E$ of the
locus
$$\{(u, v)\in E\times E-\Delta:\exists l\in G^1_3(E) \mbox{ such that }
a_1^l(q)=3,\ a_1^l(u)\geq 2, \ a_1^l(v)\geq 2\}$$ and denote by
$F_i$ the (numerical class of the) fibre of the projection $\pi_i:
E\times E\rightarrow E$ for $i=1, 2$. Using that $\cU\cap
\Delta=\{(u, u): u\neq q, \ 3 u\equiv 3 q\}$ (this intersection is
transversal!), it  follows that  $\cU\equiv 4(F_1+F_2)-\Delta$. If
$q\in E$ is viewed as the origin of $E$, then the isomorphism
$E\times E\ni (x, y)\mapsto (-x, y-x)\in E\times E$ shows that the
number of $l\in G^1_3(E)$ we are computing, equals the intersection
number $\cU\cdot \mathcal{V}$ on $E\times E$, where
$$\mathcal{V}:=\{(u, v)\in E\times E: 2 v+u\equiv 4 q-p\}.$$
Since $\mathcal{V} \equiv 3F_1+6F_2-2\Delta$, we  reach the stated
answer by direct calculation.

\noindent {\bf{(b)}} We specialize $[E, p, q]\in \cM_{1, 2}$ to the
stable curve $[E\cup_r T, p, q]\in \mm_{1, 2}$, where $[T, r, p,
q]\in \mm_{0, 3}$. We count admissible coverings
$[f:X\stackrel{4:1}\longrightarrow R, \tilde{p}, \tilde{q}]$, where
$\tilde{p}, \tilde{q}\in X_{\mathrm{reg}}$, $R$ is a nodal curve of
genus $0$ and there exist points $x, y\in X_{\mathrm{reg}}$ with the
property that the divisors $4 \tilde{p}, 3 x, \tilde{q}+2 y$ on $X$
all appear in distinct fibres of $f$. Moreover $[X, \tilde{p},
\tilde{q}]$ is a pointed curve stably equivalent to $[E\cup_r T, p,
q]$. There are three possibilities:

\noindent (1) $x, y\in E$. Then $f_T:T\stackrel{4:1}\rightarrow (\PP^1)_1$ is uniquely determined by the properties $f_T^*(0)=4p$ and $f_T^*(\infty)=3 r+q$, while $f_E:E \stackrel{3:1}\rightarrow (\PP^1)_2$ is such that $r$ and some point $x\in E-\{r\}$ appear as points of total ramification. In particular, $3x\equiv 3r$ on $E$, which gives $8$ choices for $x$. Each such $f_E$ has $2$ remaining points of simple ramification, say $y_1, y_2\in E$ and we take a rational curve $T'$ which we attach to $T$ at $q$ and map isomorphically onto $(\PP^1)_2$. Choose
$\tilde{q}\in T'$ with the property that $f(\tilde{q})=f_E(y_i)$ for $i\in \{1, 2\}$ and obviously $\tilde{p}=p\in T$.  This procedure produces $16=8\cdot 2$ admissible $\mathfrak g^1_4$'s.

\noindent (2) $x\in T,\ y\in E$. Now $f_T:T\stackrel{4:1}\rightarrow (\PP^1)_1$ has the
 properties $f_T^*(0)=4p, f_T^*(1)\geq 2r+q$ and $f_T^*(\infty)\geq 3x$ for some $x\in T$ (Up to isomorphism,
 there are $2$ choices for $f_T$). Then $f_E:E \stackrel{2:1}\rightarrow (\PP^1)_2$ is ramified
 at $r$ and at some point $y\in E-\{r\}$ such that $2y\equiv 2r$. This gives $3$ choices for $f_E$.
 We attach two rational curve $T'$ and $T''$ to $T$ at the points $q$ and $q'\in
 f_T^{-1}(f(q))-\{r, q\}$ respectively.
 We then map $T'$ and $T''$ isomorphically onto $(\PP^1)_2$. Finally we choose $\tilde{p}=p\in T$ and $\tilde{q}\in T'$
 uniquely determined by the condition $f_{T'}(\tilde{q})=f_E(y)$. We have produced $6=2\cdot 3$
 coverings.

 \noindent (3) $x\in E, y\in T$. Counting ramification points on $T$
 we quickly see that $\mbox{deg}(f_E)=3$ and $f_E:E\rightarrow
 (\PP^1)_2$ is such that $f_E^*(0)=3x$ and $f_E^*(\infty)=3r$, which
 gives $8$ choices for $f_E$. Moreover
 $f_T: T\stackrel{4:1}\rightarrow (\PP^1)_1$ must satisfy the
 properties $f_T^*(0)=4p, \ f_T^*(1)\geq q+2y$ and $f_T^*(\infty)= 3r+r'$
 for some $r'\in T$. If $[T, p, q, r]=[\PP^1, 0, 1, \infty]\in
 \mm_{0, 3}$, then $$f_T(t)=\frac{t^4}{t-r'}, \ \mbox{ where }\
 r'\in \bigl\{
 \frac{1+\sqrt{-2}}{4}, \frac{1-\sqrt{-2}}{4}\bigr\}.$$
 Thus we obtain another $16=8\cdot 2$ admissible $\mathfrak g^1_4$'s in this case.
Adding (1), (2) and (3), we found $38=16+6+16$ admissible coverings
$\mathfrak g^1_4$
 on $E\cup_r T$ and
 this finishes the proof.
\end{proof}

\begin{proposition}\label{enu3}
We fix a general pointed curve $[C, p]\in \cM_{3, 1}$. Then there
are $210$ pencils $l=\OO_C(2p+2x)\in G^1_4(C)$, $x\in C$, having an unspecified
triple point.
\end{proposition}
\begin{proof}
We define the map $\phi: C\times C\rightarrow \mathrm{Pic}^{1}(C)$
given by $$\phi(x, y):=\OO_C\bigl(2p+2x-3 y\bigr).$$ A standard
calculation  shows that $\phi^*(W_{1}(C))=g(g-1)\cdot 2^2\cdot
3^2=216$ (Use Poincar\'e's formula $[W_{1}(C)]=\theta^2/2$).
Set-theoretically it is clear that $\phi^*(W_1(C))\cap \Delta=\{(p,
p)\}$. A local calculation similar to \cite{Di} pg. 34-36, shows
that the intersection multiplicity at the point $(p, p)$ is equal to
$6=g(g-1)$, hence the answer to our question.
\end{proof}

\subsection{The invariant $N_2(d)$}
We have reached the final step of our calculation and we now compute $N_2(d)$.
We denote by $\aa_d$ the Hurwitz stack parameterizing admissible coverings of degree $d$
$$t:=[f:(X, p) \stackrel{d:1}\longrightarrow  R, \ q_0; p_0; \ p_1, \ldots, p_{6d-13}],$$
where $[X, p]$ is a pointed nodal curve of genus $2d-4$, $[R, q_0;
p_0: p_1, \ldots, p_{6d-13}]$ is a pointed nodal curve of genus $0$,
and $f$ is an admissible covering in the sense of \cite{HM} having a
point of triple ramification  $x\in f^{-1}(q_0)$, a point of simple
ramification $y\in X-\{p\}$ such that $f(y)=f(p)=p_0$ and points of
simple ramification in the fibres over $p_1, \ldots, p_{6d-13}$. The
symmetric group $\mathfrak S_{6d-13}$ acts on $\aa_d$ by permuting
the branch points $p_1, \ldots, p_{6d-13}$ and the stabilization map
 $$\phi: \aa_d/\mathfrak S_{6d-13}\rightarrow \mm_{2d-4, 1}, \mbox{  }\  \phi(t):=[X, p]$$
  is generically finite of degree $N_2(d)$.

We  completely describe the fibre $\phi^{-1}([C\cup_q E, p])$, where
$[C, q]\in \cM_{2d-5, 1}$ and $[E, q, p]\in \cM_{1, 2}$ are general
pointed curves. We count admissible covers $f:(X,
\tilde{p})\rightarrow R$ as above, where $[X, \tilde{p}]$ is stably
equivalent to $[C\cup _q E, p]$. Depending on the position of the
ramification points $x, y\in X$ we distinguish between the following
cases:

\noindent {{\bf{(i)}} $x\in C, y\in E$. From Brill-Noether theory,
we know that $\mbox{deg}(f_C)\in \{d-1, d\}$. If
$\mbox{deg}(f_C)=d$, then one possibility is that both $f_C$ and
$f_E$ are triply ramified at $q$. In this case  $f_C$ is induced by
one of the $e(d, 2d-5)$ linear series $l\in G^1_d(C)$ with
$l(-3q)\neq \emptyset$ and $l(-3x)\neq \emptyset$, for some $x\in
C-\{q\}$. The covering $f_E$ is of degree $3$ and it induces a
linear equivalence $3q\equiv 2y+p$ on $E$ which has $4$ solutions
$y\in E$. To obtain $X$ we attach to $C$ rational curves at  the
$d-3$ points in $f_C^{-1}(f(q))-\{q\}$. We have exhibited in this
way $4e(d, 2d-5)$ automorphism-free points in $\phi^{-1}([C\cup_q E,
p])$ which are counted with multiplicity $1$. Another possibility is
that both $f_C$ and $f_E$ are simply ramified at $q$ and the fibre
$f_C^{-1}(f(q))$ contains a second point $z\neq q$ of simple
ramification. The number of such $l\in G^1_d(C)$ has been denoted by
$N_3(d)$. Having chosen $f_C$, then $f_E:E\stackrel{2:1}\rightarrow
(\PP^1)_2$ is induced by $|\OO_E(2q)|$. Then we attach a rational
curve $T$ to $C$ at $z$, and we map $T\stackrel{2:1}\rightarrow
(\PP^1)_2$ using the linear system $|\OO_T(2q)|$ in such a way that
the remaining ramification point of $f_T$ maps to $f_E(p)$. We
produce $N_3(d)$ smooth points of $\aa_d/\mathfrak S_{6d-13}$ via
this construction. In both these cases $\tilde{p}=p\in C\cup E$.

\noindent {{\bf{(ii)}} $x, y\in C$.  Now  $\mbox{deg}(f_C)=d-1$ and
$f_C$ is induced by one of the $b(d-1, 2d-5)=e(d-1, 2d-5)$ linear
series $l\in G^1_{d-1}(C)$ with $l(-3x)\neq \emptyset$ for some
$x\in C-\{p\}$. Moreover, $f_C(q)$ is not a branch point of $f_C$
which implies that $\mbox{deg}(f_E)=2$ and that $f_E$ is induced by
$|\OO_E(p+q)|$. Obviously, $f_C$ and $f_E$ map to different
components of $R$. To obtain the source $(X, \tilde{p})$ of our
covering, we first attach $d-2$ rational curves to $C$ at all the
points in $f_C^{-1}(f(q))-\{q\}$ and map these curves $1:1$ onto
$f_E(E)$. Then we attach a curve $T'\cong \PP^1$, this time to $E$
at the point $q$ and map $T'$ isomorphically onto $f_C(C)$. The
point $\tilde{q}\in X$ lies on the tail $T'$ and is characterized by
the property $f_{T'}(\tilde{p})=f_C(y)$, where $y\in C$ is one of
the $6d-16$ simple ramification points of $l$. This procedure
produces $(6d-16)b(d-1, 2d-5)$ admissible coverings in
$\phi^{-1}([C\cup_q E, p])$.

\noindent {\bf{(iii)}} $x\in E, y\in E$. If $\mbox{deg}(f_C)=d$,
then $\mbox{deg}(f_E)\geq 4$ and $f_C$ is given by one of the $a(d,
2d-5)$ linear series $l\in G^1_d(C)$ such that $l(-4q)\neq
\emptyset$. Then $f_E:E\stackrel{4:1}\rightarrow \PP^1$ has the
properties that (up to an automorphism of the base) $f_E^*(0)=4q,\
f_E^*(1)\geq p+2y$ and $f^*(\infty)\geq 3x$, for some points $x,
y\in E-\{p, q\}$. The number of such $\mathfrak g^1_4$'s has been
computed in Proposition \ref{elliptic} (b) and it is equal to $38$.
Therefore this case produces $38a(d, 2d-5)$ coverings. If on the
contrary, $\mbox{deg}(f_C)=d-1$, then $f_C$ is induced by one of the
$a(d-1, 2d-5)$ linear series $l\in G^1_{d-1}(C)$ such that
$l(-2q)\neq \emptyset$,
 while $f_E:E\stackrel{3:1}\rightarrow \PP^1$ is such that (up to an automorphism of the
 base)\
 $f_E^*(0)\geq 2q, \ f_E^*(1)=p+2y,\ f_E^*(\infty)=3x$ for some $x, y\in E-\{p, q\}$.
 After making these choices, we attach $d-3$ rational curves to $C$ at the point
  $\{q'\}=f_C^{-1}(f(q))-\{q\}$
 and we map these isomorphically onto $f_E(E)$. Furthermore, we attach a rational curve $T'$ to $E$ at
 the point  $\{q'\}=f_E^{-1}(f(q))-\{q\}$ and map $T'$ isomorphically onto $f_C(C)$. Using Proposition
 \ref{elliptic} (a), we obtain $11a(d-1, 2d-5)$ admissible coverings. Altogether part (iii)
provides $38a(d-1, 2d-5)+11a(d-1, 2d-5)$ points in $\aa_d/\mathfrak
S_{6d-13}$.

\noindent {{\bf{(iv)}} $x\in E, y\in C$. In this case, since $p$ and
$y$ lie in different components, we know that we have to ``blow-up''
the point $p$ and insert a rational curve which is mapped to the
component $f_C(C)$ of $R$. Thus $\mbox{deg}(f_{C})\leq d-1$, and by
Brill-Noether theory it follows that $\mbox{deg}(f_C)=d-1$.
Precisely, $f_C$ is induced by one of the $a(d-1, 2d-5)$ linear
series $l\in G^1_{d-1}(C)$ such that $l(-2q)\neq \emptyset$.
Furthermore, $f_E:E\stackrel{3:1}\rightarrow \PP^1$ can be chosen
such that $f_E^*(0)=p+2q$ and $f_E^*(\infty)=3x$ for some $x\in E$.
This gives the linear equivalence $3x\equiv p+2q$ on $E$ which has
$9$ solutions. We attach $d-3$ rational curves at the points in
$f_C^{-1}(f(q))-\{q\}$ and map these $1:1$ onto $f_E(E)$. Finally,
we attach a rational curve $T'$ to $E$ at the point $p$ and map $T'$
such that $f(T')=f(C)$. We pick $\tilde{p}\in T'$ with the property
that $f_{T'}(\tilde{p})=f_C(y)$, where $y\in C$ is one of the
$6d-15$ ramification points of $f_C$. We have obtained
$9(6d-15)a(d-1, 2d-5)$ admissible coverings in this way.

We have completely described $\phi^{-1}([C\cup_q E, p])$ and it is
easy to check that all these coverings have no automorphisms, hence
they give rise to smooth points in $\aa_d$ and that the map $\phi$
is unramified at each of these points. Thus
$$N_2(d)=\mbox{deg}(\phi)
=4e(d, 2d-5) + (6d-16)b(d-1, 2d-5)+38a(d, 2d-5)+$$
$$+11 a(d-1, 2d-5) +9(6d-15)\ a(d-1, 2d-5)+ N_3(d).$$
For $d=4$, we know that $N_3(4)=210$ (cf. Proposition \ref{enu3}),
which determines $N_2(4)$ and the class $[\dd_3]$. We record these
results:
\begin{theorem}\label{d3}
The locus $\mathcal{D}_3$ of pointed curves $[C, p]\in \cM_{2, 1}$
with a pencil $l\in G^1_4(C)$ totally ramified at $p$ and having two
points of triple ramification, is a divisor on $\cM_{2, 1}$. The
class of its compactification in $\mm_{2, 1}$ is given by the
formula:
$$\dd_3 \equiv 640\psi-860\lambda+72\delta_0\in \mathrm{Pic}(\mm_{2, 1}).$$
\end{theorem}
\begin{theorem}\label{n2}
For a general pointed curve $[C, p]\in \cM_{2d-4, 1}$ the number of
pencils $L\in W^1_d(C)$ satisfying the conditions
$$h^0(L\otimes \OO_C(-3x))\geq 1\  \mbox{ and }\ h^0(L\otimes
\OO_C(-p-2y))\geq 1$$ for some points $x, y\in C-\{p\}$, is equal to
$$N_2(d)=\frac{6(40d^2-179d+212)\ (2d-4)!}{d!\ (d-3)!}\ .$$
\end{theorem}
\begin{remark} As a check, for $d=3$, the number $N_2(3)$ computes
the number of pairs $(x, y)\in C\times C$ such that $p\neq x\neq
y\neq p$ and $3x\equiv p+2y$. This number is equal to $r(3, 2)=70$
which matches Theorem \ref{n2}.
\end{remark}
\begin{theorem}\label{nn3}
We fix an integer $d\geq 4$. For a general pointed curve $[C, p]\in
\cM_{2d-5, 1}$, the number of pencils $L\in W^1_d(C)$ satisfying the
conditions
$$h^0(L\otimes \OO_C(-3x))\geq 1\ \mbox{ and }\ h^0(L\otimes
\OO_C(-2p-2y))\geq 1$$ for some points $x, y\in C-\{p\}$, is equal
to
$$N_3(d)=\frac{84(d-3)(2d^2-10d+13)\ (2d-4)!}{d!\ (d-2)!}\ .$$
\end{theorem}
\begin{remark} For $d=4$, Theorem \ref{nn3} specializes to
Proposition \ref{enu3} and we find again that $N_3(4)=210$.
\end{remark}

\end{document}